\documentclass{amsart}

\usepackage{latexsym}
\usepackage{oldgerm}
\usepackage{mathrsfs}
\usepackage{amsfonts}
\usepackage{amsmath}
\usepackage{amssymb}
\def\b{\mathbb }
\def\scr{\mathscr}
\def\phi{\varphi }

\newtheorem{theorem}{Theorem}[section]
\newtheorem{lemma}[theorem]{Lemma}
\newtheorem{proposition}[theorem]{Proposition}

\newtheorem{corollary}[theorem]{Corollary}
\theoremstyle{definition}
\newtheorem{definition}[theorem]{Definition}

\newtheorem{examples}[theorem]{Examples}

\theoremstyle{remark}
\newtheorem{remark}[theorem]{Remark}
\newtheorem{remarks}[theorem]{Remarks}

\numberwithin{equation}{section}

\begin{document}

\title{A positive radial product formula for the Dunkl kernel}

\author{Margit R\"osler}
\address{Mathematisches Institut, Universit\"at G\"ottingen\\ 
Bunsenstr. 3--5, D-37073 G\"ottingen\\
Germany}
\email{roesler@uni-math.gwdg.de}


\date{October 8, 2002}


\begin{abstract}
It is an open conjecture 
that generalized Bessel functions associated with root systems have
a positive product formula for non-negative multiplicity parameters of 
the associated Dunkl operators. 
In this paper, a partial result towards this conjecture is proven, namely a positive radial
product formula for the non-symmetric counterpart of 
the generalized Bessel function, the Dunkl kernel. Radial hereby means that one of the factors
in the product formula is replaced by its mean over a sphere. 
The key to this product formula is a positivity result for the Dunkl-type spherical mean operator.
It
 can also be interpreted in the sense that the Dunkl-type generalized 
translation of radial functions
is positivity-preserving. As an application, we construct 
Dunkl-type homogeneous  Markov processes associated with radial probability distributions.

\end{abstract}

\maketitle


\section{Introduction}

Along with addition formulas, product formulas have always been a 
challenging topic in the area of one-variable special functions. Typically, 
positive product formulas are obvious for particular parameters 
from a group theoretical
background (when e.g. the  functions under consideration have an interpretation as the
spherical functions of a Gelfand pair), but it is difficult to obtain
a generalization to larger classes of parameter values. For special functions in several
variables such questions 
seem to be even more intricate. One open conjecture in this direction
concerns  positive product formulas for multivariable Bessel functions associated with root
systems, provided all multiplicity parameters are nonnegative.
 In the present paper, we 
prove a partial result towards a positive product formula for multivariable Bessel functions
on a Weyl chamber, namely a positive radial
product formula for the non-symmetric counterpart of 
the generalized Bessel function, the Dunkl kernel. This kernel is the analogue of the usual
exponential function in the theory of rational Dunkl operators as 
developed in \cite{D1} - \cite{D4}.
To become more precise, let us briefly introduce our setting.
 Let $R$ be a (reduced, not necessary crystallographic)
 root system in $\b R^N,$ equipped with the standard Euclidean inner product $\langle\,.\,,\,.\,\rangle,$  i.e. R$\subset \b R^N\setminus\{0\}$ is finite with $R\cap \b R\alpha = \{\pm\alpha\}$  and $\sigma_\alpha(R) = R$ for all $\alpha\in R$, where $\sigma_\alpha$ denotes
 the reflection in the hyperplane perpendicular to $\alpha$. We assume 
$R$ to be normalized such that $\langle\alpha,\alpha\rangle = 2$ for all 
$\alpha\in R$. This simplifies formulas, but is no loss of generality 
for our purposes. 
Further let $G$ denote the finite reflection group generated by $\{\sigma_\alpha\,,\,\alpha\in R\}$
and let
$k: R\to \b C$  be a fixed
multiplicity function on $R,$  i.e. a function
which is constant on the orbits under the action of $G$. We shall always assume that $k$ is non-negative, i.e. $k(\alpha) \geq 0$ for all $\alpha\in R$. 
The  first-order rational Dunkl operators attached to
$G$ and $k$ are defined by
\begin{equation}\label{(1.10a)}
 T_\xi(k) f(x)  \,=\, \partial_\xi f(x)  + 
\sum_{\alpha\in R_+} k(\alpha) \langle \alpha,\xi\rangle 
\frac{f(x) - f(\sigma_\alpha x)}{\langle\alpha,x\rangle},\quad x,\,\xi\in \b R^N.
\end{equation}
Here $\partial_\xi$ denotes the derivative in direction $\xi$ and $R_+$ is some fixed positive subsystem  of $R$. The definition is independent of the special choice of $R_+$, thanks to the $G$-invariance of $k$. As first shown in \cite{D1}, the 
$T_\xi(k), \,\xi\in \b R^N$ 
 generate a commutative algebra of differential-reflection operators.  This is the foundation for  rich analytic structures related
with them.
In particular, there  exists a   counterpart of the usual exponential function, called the Dunkl kernel, and an analogue of the Euclidean Fourier transform
with respect to this kernel. The Dunkl kernel $E_k$ is  holomorphic 
on $\b C^N\times \b C^N$ and symmetric in its arguments. 
  Similar to the situation for
spherical functions
on a Riemannian symmetric space, the function $E_k(\,.\,,y)$ with fixed 
spectral parameter $y\in \b C^N$ may  be 
characterized as the unique analytic solution of the
 joint eigenvalue problem 
\begin{equation}\label{(1.104)}
 T_\xi(k) f\,=\, \langle\xi,y\rangle f \quad\text{for all }\,\xi \in \b C^N,\,\,\,
f(0) = 1,
\end{equation} c.f. \cite{O}.
Apart from the trivial case $k=0$, where $E_k(x,y) = e^{\langle x,y\rangle}$, the
kernel $E_k$  is explicitly known in very few cases only. These include the
rank-one case (see Section \ref{rank-one}.) as well as the symmetric group $G=S_3$ 
(\cite{D3}).
The reflection invariant counterpart of $E_k$ is the so-called generalized Bessel function
\[ J_k(x,y) = \frac{1}{|G|} \sum_{g\in G} E_k(gx,y),\]
which has been first studied in  \cite{O}.  It is $G$-invariant in both arguments 
and naturally considered on the Weyl chambers of $G$. If the rank of $R$ is one, then $J_k$ 
 coincides with a usual Bessel function.
In particular cases,
all being related with Weyl groups and certains
half-integer multiplicity parameters,
generalized Bessel functions can be given an interpretation as the  spherical
functions of a Euclidean-type symmetric space. We shall have a closer look on these
 examples in the appendix, Section \ref{euclidean}.
In  these cases, and for all non-negative multiplicity-parameters in the
rank-one case, the generalized Bessel functions $J_k$ 
 have a positive product formula of
the form
\begin{equation}\label{(1.1f)}
 J_k(x,z)J_k(y,z)\,=\, \int_{\overline C} J_k(\xi,z)
         d\nu_{x,y}^k(\xi) \quad\text{for all }z\in \b C^N.\end{equation}
Here $\overline C$ is the topological closure of the
Weyl chamber 
\[ C = \{ x\in \b R^N: \langle\alpha,x\rangle > 0 \quad\text{for all }\, 
\alpha\in R_+\},\]
and the
  $\nu_{x,y}^k$ are compactly supported probability measures on 
 $\overline C$.
It is conjectured 
 that  a positive product formula of the form \eqref{(1.1f)}
should in fact be valid for 
arbitrary reflection groups and non-negative multiplicity parameters. 
In that case, one would obtain a positivity-preserving convolution of
regular bounded Borel measures on the chamber $\overline C$ when defining the convolution of 
point measures according to 
$\delta_x\circ_k\delta_y:= \nu_{x,y}^k$.  In the special cases mentioned above,
this convolution induces the structure of a commutative hypergroup on
 $\overline C$,
 and we conjecture that this should be 
true in general (for non-negative multiplicities). Roughly speaking, a hypergroup consists of 
a locally compact Hausdorff space $H$  together with a positivity-preserving convolution on the space of regular bounded Borel measures on $H$, which allows one to carry
over the harmonic analysis on a locally compact group to a large extent. 
The hypergroup setting in particular includes double coset spaces $G//H$ with $H$ a compact subgroup of a locally compact group $G$. For an introduction to the subject, the reader is referred to \cite{BH} or \cite{J}.

As already indicated, the generalized exponential function $E_k$ is of particular interest as it
gives rise to an integral transform on
$\b R^N$, commonly called the Dunkl transform. On suitable function spaces, 
this transform
maps  Dunkl operators to multiplication operators. The Dunkl transform
 provides a natural generalization of the usual Euclidean 
Fourier transform, to which it reduces in case $k=0$. 
 It is another  open question 
whether the Dunkl transform  admits an interpretation as the Gelfand transform
of a suitable $L^1$-convolution algebra  on $\b R^N$
for arbitrary parameters $k\geq 0$. Again this is true in the rank-one case.
Here
the Dunkl kernel satisfies a
 product formula which leads to a convolution structure on the
entire real line, providing  a  natural extension of the usual  group
structure, see \cite{R0}. In contrast to a hypergroup convolution, this convolution is not positivity-preserving if $k>0$.
We conjecture that an analogous statement is true for arbitrary rank and all
multiplicities $k\geq 0$. In particular,  the 
Dunkl kernel should have a
 product formula of the form
\begin{equation}\label{(1.1g)}
 E_k(x,z)E_k(y,z)\,=\, \int_{\b R^N} E_k(\xi,z)
         d\mu_{x,y}^k(\xi) \quad\text{for all }z\in \b C^N\,,
\end{equation}
where the measures $\mu_{x,y}^k$ are signed Borel measures  on $\b R^N$ which 
are uniformly bounded with respect to total variation norm.

This paper presents a result towards both the stated conjectures, 
namely a positive 
``radial'' product formula for $E_k$.
More precisely, we shall prove that for each $x\in \b R^N$ and $t\geq 0$ there
exists a unique compactly supported probability measure $\sigma_{x,t}^k$ on $\b R^N$ such that
\begin{equation}\label{(1.1c)}
E_k(ix,z) j_{\lambda}(t|z|)\,=\,\int_{\b R^N} E_k(i\xi,z) d\sigma_{x,t}^k(\xi)
\quad\text{for all }\, z\in \b R^N.
\end{equation}
Hereby the index of the Bessel function $j_{\lambda}$ is 
given by $\lambda = \gamma +N/2-1$, with 
\[\gamma = \sum_{\alpha\in R_+} k(\alpha) \geq 0.\]
Of course, the kernel $E_k$ in \eqref{(1.1c)} may be equally replaced by the generalized Bessel function $J_k$.
This product formula will be obtained from a study of  
the Dunkl-type spherical mean operator, as first introduced in \cite{MT}.
In analogy to the classical case, the spherical mean operator $f\mapsto M_f$ 
is defined for $f\in C^\infty(\b R^N)$ by
\[ M_f(x,t) = 
\frac{1}{d_k} \int_{S^{N-1}} 
      f(x*_kty) w_k(y) d\sigma(y) \quad(x\in \b R^N, \, t\geq 0).\]
Here $d_k$ is a normalization constant, 
$w_k$ is the $G$-invariant  weight function
\begin{equation}\label{(1.103)} w_k(x)\,=\, \prod_{\alpha\in R_+} |\langle
\alpha,x\rangle|^{2k(\alpha)} \quad(x\in \b R^N)
\end{equation}
and $f(x*_ky)$ denotes the Dunkl-type generalized translation, which coincides with the
usual group translation on $\b R^N$ in case $k=0$ and satisfies
\[ E_k(x*_k y,z) = E_K(x,z) E_k(y,z) 
\quad\text{for all }\,x,y \in \b R^N, \, z\in \b C^N.\]
The key result of the present paper states that for $k\geq 0$, the
associated spherical mean operator is positivity-preserving. This implies
the existence of compactly supported probability measures $\sigma_{x,t}^k$
on $\b R^N$ which represent this operator in the sense that
\[ M_f(x,t)\,=\, \int_{\b R^N} fd\sigma_{x,t}^k \quad\text{for all }\,
f\in C^\infty(\b R^N).\]
When specializing to $f(x) = E_k(ix,z)$ with $z\in \b R^N$, one obtains 
\eqref{(1.1c)}. 
Similar to the classical case $k=0,$  $M_f$ satisfies a second order
differential-reflection equation of Darboux-type.
A study of the domain
of dependence for this equation allows one to deduce further information
on the support of the representing measures $\sigma_{x,t}^k$. In fact, 
the support of $\sigma_{x,t}^k$ turns out to be contained in the union of closed
balls with
radius $t$ around the points $gx, \, g\in G$. In contrast to the classical 
case, where the support reduces to the sphere $\{\xi\in \b R^N: |\xi-x| = t\},$
now  full balls as well as the complete 
 $G$-orbit of $x$ has to be taken into account. This is
due to the reflection parts 
in the generalized Darboux equation for $M_f$.  
 These results are
contained in Theorem \ref{T:Main}.
A slightly weaker variant of 
Theorem \ref{T:Main} is given in Theorem \ref{rad}. Here it is shown that 
the Dunkl-type generalized translation  (more precisely, the mapping
$f\mapsto f(x*_ky)$) is 
positivity-preserving when restricted to radial $f$. In Section \ref{S:Markov}, the 
obtained results are applied to construct ``radial'' 
semigroups $(P_t)_{t\geq 0}$ of Markov kernels
on $\b R^N$, which are translation invariant in the generalized sense 
of \cite{RV}. Hereby $P_t$ is obtained by a Dunkl-type translation of a probability 
measure $\mu_t$ of the form $d\mu_t(x) = w_k(x) d\mu_t^\prime(x),$ where  
$\mu_t^\prime$ is rotation invariant.

\section{Preliminaries}
\label{sec:1}

In this introductory part we give an account on results from Dunkl theory which will be relevant for the sequel. These concern in particular the Dunkl kernel, the  Dunkl transform and 
generalized translations. The series expansion for the Dunkl kernel according derived in Section 
\ref{S:Series} is of some interest in its own and might be known, but it is, 
to the author's knowledge, nowhere written down.
We also include a discussion of the rank-one case as a motivating 
example. For a further background in Dunkl theory,  the reader is referred to \cite{D1} -- \cite{D4}, \cite{dJ1}, \cite{R2}, \cite{T3}, and \cite{DX}. Concerning root systems and reflection groups, see \cite{Hu}. 

Throughout the paper,  $\langle\,.\,,\,.\,\rangle$  denotes the standard Euclidean scalar product in $\b R^N$ as well as its bilinear extension to $\b C^N\times \b C^N$. For $x\in \b R^N$, we write
 $|x|= \sqrt{\langle x,x\rangle}.$ Further, 
 $\b Z_+:= \{0,1,2, \ldots\}$, and $\b R_+:=
[0,\infty)$.
We denote by  $C^\infty(\b R^N)$  the space of infinitely often differentiable functions on $\b R^N$  and by
 $\mathscr S(\b R^N)$ the  Schwartz space of rapidly decreasing functions, both 
equipped with the usual
Fr\'echet space topologies. For a locally compact Hausdorff space $X,$ $C_b(X)$ 
denotes
the space of continuous, bounded functions on $X$. Further, $M_b(X),\, M_b^+(X)$ 
and $M^1(X)$ stand for the spaces of regular bounded complex Borel measures on $X,$ those which are non-negative, and those which are probability measures respectively. The 
$\sigma(M_b(X), C_b(X))$-topology is referred to as  the weak topolgy on $M_b(X),$ and
the $\sigma$-algebra of Borel sets in $X$ is denoted by $\mathcal B(X)$.

\subsection{Basics from Dunkl theory}

Let $G$ be a finite reflection group on $\b R^N$ with root system $R$, and
fix a positive subsystem $R_+$ of $R$ as well as a  
non-negative multiplicity 
function $k$. The associated Dunkl operators $T_{\xi}(k)$, defined according
to \eqref{(1.10a)},
share many properties with usual
partial derivatives. 
In particular,  if $f\in C^k(\b R^N),$ then $T_\xi(k)f\in C^{k-1}(\b R^N),$  the $T_\xi(k)$ are homogeneous of degree $-1$ on polynomials and they
 leave the Schwartz space $\mathscr S(\b R^N)$ invariant. Moreover, $T_\xi(k)$ is
  $G$-equivariant:
\[ g\circ T_\xi(k)\circ g^{-1} = T_{g\xi}(k)  \quad(g\in G).\]
The counterpart of the usual Laplacian is the Dunkl Laplacian, defined by
\[ \Delta_k := \sum_{i=1}^N T_{\xi_i}(k)^2,\]
where $\{\xi_i\,,i=1,\ldots,N\}$ is an arbitrary orthonormal basis of
$(\b R^N, \langle\,.\,,\,.\,\rangle),$ c.f. \cite{D1}.
It is given explicitly by 
\[ \Delta_k f(x) \,=\, L_k f(x) - 2\sum_{\alpha\in R_+} k(\alpha) \frac{f(x) - 
f(\sigma_\alpha x)}
{\langle\alpha,x\rangle^2},\]
with the singular elliptic operator
\begin{equation}\label{(2.108)}
 L_kf(x):=\,  \Delta f(x)  + 2\sum_{\alpha\in R_+}
k(\alpha)\frac{\langle\nabla f(x),\alpha\rangle}
{\langle\alpha ,x\rangle}.
\end{equation}

According to \cite{D2}, there exists 
    a unique,  degree of  homogeneity preserving linear 
isomorphism $V_k$ on polynomials such that
\[T_\xi(k)V_k\,=\, V_k\partial_\xi \quad\text{for all }\, 
\xi \in \b R^N \,\,\text{ and }\, 
V_k(1) = 1. \]
It is shown  in  \cite{R2} that  $V_k$ has a Laplace-type 
representation of the form 
\begin{equation}\label{(1.40a)}
V_k f(x)\,=\, \int_{\b R^N} f(\xi)d\mu_x^k(\xi)
\end{equation}
with a unique probability measure $\mu_x^k\in M^1(\b R^N)$  whose support 
is contained in 
\[C(x) := \text{co}\{gx,g\in G\},\]
 the convex hull of the $G$-orbit of $x$ in $\b R^N$. 
By means of formula \eqref{(1.40a)}, $V_k$ may be extended  to various larger function spaces including $C^\infty(\b R^N).$ We denote this extension  by $V_k$ again.
In fact, $V_k$ establishes a homeomorphism of 
$C^\infty(\b R^N)$, see Theorem 4.6 of \cite{dJ2} or \cite{T3}.

 The Dunkl kernel asociated with $G$ and $k$ is defined by
\begin{equation}\label{(2.77)}
 E_k(x,y)=  V_k(e^{\langle\, .\,, y\rangle})(x)\,=\, \int_{\b R^N} 
e^{\langle \xi, y\rangle}d\mu_x^k(\xi) \quad (x\in \b R^N, \, y\in \b C^N).
\end{equation}
For fixed $y$, $E_k(\,.\,,y)$ is the unique real-analytic solution of 
\eqref{(1.104)}, c.f. \cite{O}.
 The kernel $E_k$ 
 is symmetric in its arguments and has a unique
holomorphic extension to $\b C^N\times \b C^N.$ Moreover, 
\begin{equation}\label{(2.102)}
 E_k(\lambda z,w) = E_k(z,\lambda w) \quad\text{and }\, E_k(gz,gw) = E_k(z,w)
\end{equation}
for all $z,w\in \b C^N, \,\lambda\in \b C$ and $g\in G$. 
 Let $w_k$ denote the $G$-invariant weight function \eqref{(1.103)}. 
The associated Dunkl transform  on $L^1(\b R^N, w_k)$ is then defined by
\[ \widehat f^{\,k}(\xi):=  c_k^{-1}\int_{\b R^N} f(x) E_k(-i\xi,x) 
\,w_k(x)dx \>\>(\xi\in \b R^N).\]
Here $c_k$ denotes the Mehta-type constant
\[ c_k:= \int_{\b R^N} e^{-|x|^2/2} w_k(x)dx.\]
We shall also consider the Dunkl transform on the measure space $M_b(\b R^N),$ 
\[ \widehat \mu^{\,k}(\xi):= \int_{\b R^N} E_k(-i\xi,x) d\mu(x), \quad (\xi\in \b R^N).\]
Many properties of the Euclidean Fourier transform carry over to
the Dunkl transform. The results listed below can be found in 
\cite{D4}, \cite{dJ1} and \cite{RV}:

\begin{proposition}\label{Dunkltrafo}
\begin{enumerate}
\item[\rm{(1)}] The Dunkl transform $f\mapsto \widehat f^{\,k}$ 
        is a homeomorphism of 
$\mathscr S(\b R^N)$.
\item[\rm{(2)}] ($L^1$-inversion) If $f\in L^1(\b R^N, w_k)$ with $\widehat f^{\,k}
\in L^1(\b R^N, w_k),$ then $\,f = (\widehat f^{\,k}\,)^{\vee k}\,$ a.e.
Its inverse is given by $\, f^{\vee k}(\xi) := \widehat f^{\,k}(-\xi)$. 
\item[\rm{(3)}] (Plancherel's Theorem) The Dunkl transform on $\mathscr S(\b R^N)$ extends 
uniquely to an isometric isomorphism of $L^2(\b R^N, w_k).$ 
\item[\rm{(4)}] The Dunkl transform is injective on $M_b(\b R^N).$ 
\item[\rm{(5)}] (L\'evy's continuity theorem) Let
 $(\mu_n)_{n\in \b N}\subset M_b^+(\b R^N)$ such that $(\widehat\mu^{\,k}_n)_{n\in \b N}$
converges pointwise to a function $\phi:\b R^N\to \b C$  which is continuous at $0$. 
Then there exists a unique measure $\mu\in M_b^+(\b R^N)$ with $\widehat \mu^{\,k} =\phi,$ and
$(\mu_n)_{n\in \b N}$ tends weakly to $\mu$.   
\end{enumerate}
\end{proposition}

In  \cite{T3}, a generalized translation  on  
$C^\infty(\b R^N)$ is defined by
\begin{equation}\label{(1.5)}
 \tau_yf(x)  := V_k^x V_k^y (V_k^{-1}f)(x+y),\quad x, y\in \b R^N.
\end{equation} 
Hereby the uppercase index denotes the relevant variable. 
Notice that
\[\tau_0 f = f,\, T_\xi(k)\tau_y f = 
     \, \tau_y T_\xi(k)f \quad\text{and }\, \tau_y f(x) = \tau_x f(y)\quad\text{for all }\, 
x,y\in \b R^N.\]
We shall frequently use the more suggestive notation 
\[ f(x*_k y):= \tau_y f(x).\]
For $k=0$, one just obtains the usual group translation on $\b R^N$:
$\,f(x*_0y) = f(x+y).$ 
It is also immediate from the definition that
\begin{equation}\label{(1.5b)}
E_k(x*_ky,z) = E_k(x,z) E_k(y,z) \quad\text{for all }\,z\in \b C^N.
\end{equation}
We collect some further properties of this translation which will be used later on;
for the proofs, the reader is referred to \cite{T3}.

\begin{lemma} \label{L:Translation}
\begin{enumerate}
\item[\rm{(1)}] For fixed $y\in \b R^N, \,\,\tau_y$ is a continuous linear mapping
from $C^\infty(\b R^N)$ into $C^\infty(\b R^N)$. 
\item[\rm{(2)}] For fixed $x,y\in \b R^N,$ the mapping $f\mapsto f(x*_ky)$ defines
 a compactly supported 
distribution. Its support is contained in the ball
$\{\xi\in \b R^N: |\xi|\leq |x|+|y|\}.$
\item[\rm{(3)}] If $f\in \mathscr S(\b R^N),$ then also
$\,\displaystyle  \tau_y f\,\in \mathscr S(\b R^N), $ and  $\,\displaystyle\, 
 (\tau_y f)^{\wedge k}(\xi)\,=\,
  E_k(iy,\xi)\, \widehat f^{\,k}(\xi).$ Moreover,
\[\tau_y f(x)=\, \frac{1}{c_k} \int_{\b R^N}  \widehat f^{\,k}(\xi) 
\,E_k(ix,\xi) E_k(iy, \xi)\, w_k(\xi) d\xi .\]
\item[\rm{(4)}] If $f,\, g\in \mathscr S(\b R^N)$ and $x\in \b R^N$, then
\[ \int_{\b R^N} f(x*_ky)g(y)w_k(y)dy\,=\, 
\int_{\b R^N} f(y)g(-x\,*_k y)w_k(y)dy\,.\]
\end{enumerate}
\end{lemma}

\begin{remarks}
Property (3) reveals that on $\mathscr S(\b R^N),$ the translation 
\eqref{(1.5)} coincides with the version previously introduced in 
\cite{R1}. In \cite{T3}, part (4) is shown only for compactly supported 
test functions;
a simple density argument gives the result for Schwartz functions. Alternatively,
(4) follows immediatley from (3) and the Plancherel theorem for the Dunkl transform.
\end{remarks}


\subsection{Expansion of $E_k$ in terms of $k$-spherical harmonics}\label{S:Series}

In this section, we derive a series representation for the Dunkl 
kernel 
in terms of generalized spherical (``$k$-spherical'') harmonics, which will be 
employed in the positivity proof for the Dunkl-type spherical mean operator. 
For a background in $k$-spherical harmonics, 
the reader may consult the recent monograph \cite{DX}. Throughout this section it is assumed that $N\geq 2$. 
The
 space of $k$-spherical harmonics of degree $n\geq 0$ is  defined by 
\[ \mathcal H_n^k = 
\text{ker} \Delta_k \cap \mathcal P_n^N\,,\]
where $\Delta_k$ is the Dunkl Laplacian and 
$\mathcal P_n^N$ denotes the space of homogeneous polynomials of 
degree $n$ on $\b R^N$. The space $\mathcal H_n^k$ has a reproducing kernel 
$P_n^k(\,.\,,\,.\,),$ which  is defined by 
 the property
\[ f(x) \,=\, d_k^{-1}\int_{S^{N-1}}f(y) P_n^k(x,y) w_k(y)d\sigma(y) 
 \quad\text{ for all }\,f\in \mathcal H_n^k \quad\text{and }\, |x|<1.\]
 Here $S^{N-1} = \{x\in \b R^N: |x|=1\}$ is the unit sphere in $\b R^N, \,\,
d\sigma$ denotes the 
Lebesgue surface measure and 
\[d_k = \int_{S^{N-1}} w_k(x)d\sigma(x)\, =\, \frac{c_k}{2^\lambda\Gamma(\lambda+1)}, \] 
with 
\[\lambda = \gamma +N/2-1 \geq 0.\]

\noindent
Suppose that  
 $\{Y_{n,j}:\, j=1,\ldots, d(n,N)\},\,$  is a real-coefficient orthonormal basis  of $\mathcal H_n^k$ 
in
$L^2(S^{N-1}, d_k^{-1}w_kd\sigma)$. 
In terms of this basis,  $P_n^k$ is given by
\[ P_n^k(x,y)\,=\, \sum_{j=1}^{d(n,N)} Y_{n,j}(x) Y_{n,j}(y)\,.\]
 If $x,y\in S^{N-1}$ then, according to Theorem 3.2 of \cite{X1},  
the kernel $P_n^k(x,y)$ can be written as 
\begin{equation}\label{(1.1b)}
 P_n^k(x,y) \,=\, \frac{(n+\lambda)(2\lambda)_n}{\lambda\cdot n!}\, 
 V_k\widetilde C_n^{\lambda}(\langle x,\,.\,\rangle)(y),
\end{equation}
where the $\widetilde C_n^{\lambda} \,(n\in \b Z_+)\,$ 
are the (renormalized)
Gegenbauer polynomials   
\begin{align}
\widetilde C_n^\lambda(x)\,=&\, \frac{(-1)^n}{2^n(\lambda+1/2)_n}
(1-x^2)^{1/2-\lambda} \frac{d^n}{dx^n}(1-x^2)^{n+\lambda-1/2}\notag\\ 
\,=&\, _2F_1\bigl(-n,n+2\lambda,\lambda+1/2; \frac{1-x}{2}\bigr).
\end{align}
Notice that in terms of this normalization, formula \eqref{(1.1b)}
remains true in the limiting case $\lambda = 0$.

\begin{proposition} Let $N\geq 2$. Then for all $x,y\in \b R^N$ the Dunkl kernel 
$E_k(ix,y)$ admits the  representation
\begin{equation}\label{(1.2)}
E_k(ix,y)\,=\, 
\sum_{n=0}^\infty \frac{\Gamma(\lambda +1)}{2^n\Gamma(n+\lambda+1)} 
j_{n+\lambda}(|x||y|)\,P_n^k(ix,y),
\end{equation}
the convergence of the series being uniform on compact subsets of 
$\b R^N\times \b R^N$.
\end{proposition} 

\begin{proof} 
By Gegenbauer's degenerate form of the addition theorem
 for 
Bessel functions (\cite{Wa}, p. 368), we have
\begin{equation}\label{(1.3)} 
e^{irt}\,=\,\sum_{n=0}^\infty \bigl(\frac{ir}{2}\bigr)^n 
\frac{(2\lambda)_n}{(\lambda)_n\,n!}\, j_{n+\lambda}(r) 
\widetilde C_n^\lambda(t) 
\quad \text{for all } r\in \b R,\, t\in [-1,1]
\end{equation}
(With the obvious extension to the case $\lambda = 0$).
The series converges uniformly on every compact subset of 
$\b R \times [-1,1]$. This is easily seen from the asymptotic 
behavior of the Gamma function and the estimates 
\[|j_{n+\lambda}(r)| \leq 1,  \quad \quad
|\widetilde C_n^\lambda(t)|\leq \widetilde C_n^\lambda(1)\, = \,1\]
which hold within the relevant ranges of $r$ and $t.$ 
We now put $t=\langle x,y\rangle$ with $x,y\in S^{N-1}$, and apply the intertwining operator $V_k$ to both sides of \eqref{(1.3)}. This may be done termwise, the locally uniform convergence being maintained: in 
fact, according to \eqref{(2.77)} together with
\eqref{(1.3)} and \eqref{(1.1b)}, we obtain
\begin{align} E_k(irx,y)\,=&\,  V_k(e^{ir\langle x,\,.\,\rangle})(y)
\,=\,\int_{|\eta|\leq 1} e^{ir\langle x,\eta\rangle} d\mu_{y}^k(\eta) \notag\\
&= \sum_{n=0}^\infty \frac{(2\lambda)_n}{(\lambda)_n\,n!}
\bigl(\frac{ir}{2}\bigr)^n j_{n+\lambda}(r)\, V_k 
\widetilde C_n^\lambda(\langle x,\,.\,\rangle)(y) \notag\\ &=\,  \sum_{n=0}^\infty \frac{\Gamma(\lambda+1)}{2^n\Gamma(n+\lambda+1)}\, j_{n+\lambda}(r) 
P_n^k(irx,y)\,.\notag
\end{align}
As $\, |V_k \widetilde C_n^\lambda(\langle x,\,.\,\rangle)(y)|\,\leq \widetilde C_n^\lambda(1)\,= 1\,$ for all $x,y\in  S^{N-1}$, the series converges uniformly on compact subsets of 
$\b R\times S^{N-1}\times S^{N-1}$.
This implies the assertion. 
\end{proof}
 
\noindent
In view of the orthogonality of the $k$-spherical harmonics 
$Y_{n,j}$, termwise spherical integration of \eqref{(1.2)}  
leads to  the following well-known
special case of the Dunkl-type Funk-Hecke formula 
(Theorem 5.3.4 of \cite{DX}):

\begin{corollary}\label{Funk} Let $N\geq 2.$ Then for all $x\in \b R^N,$ 
\[ \frac{1}{d_k}\int_{S^{N-1}} E_k(ix,y) Y_{n,j}(y) w_k(y)d\sigma(y)\,=
\,\frac{\Gamma(\lambda +1)}{2^n\Gamma(n+\lambda+1)}j_{n+\lambda}(|x|) 
   Y_{n,j}(ix)\,.\]
\end{corollary}

\subsection{The rank-one case}\label{rank-one}

Here the reflection group is $G=\{id,\sigma\}$, acting on $\b R$ via 
$\sigma(x) = -x$. The corresponding kernel
$E_k$ with parameter $k\geq 0$ has been calculated in 
\cite{D3}. It is given by
\[ E_k(z,w)\,=\, j_{k-1/2}(izw)\,+\,\frac{zw}{2k+1}\,
j_{k+1/2}(izw) \] 
with  the normalized spherical Bessel function
\[ j_\alpha(z)\,=\,
\Gamma(\alpha+1)\cdot\sum_{n=0}^\infty
\frac{(-1)^n(z/2)^{2n}}{n!\,\Gamma(n+\alpha+1)}\,.\]
Thus $J_k(z,w) = j_{k-1/2}(izw).$ 
It is well-known (see e.g. \cite{BH}, 3.5.61) that the Bessel functions $j_\alpha$ with $\alpha\geq -1/2$  
 satisfy a  product formula of the form
\begin{equation}\label{(2.106)}
j_\alpha(xz) j_\alpha(yz)\,=\,\int_0^\infty j_\alpha(\xi z) d\nu_{x,y}^\alpha(\xi)
\quad\text{for all }\,z\in \b C.
\end{equation}
Here the  $\nu_{x,y}^\alpha$ are probability measures on $\b R_+$. 
For $x,y >0$, they are given by
\[ d\nu_{x,y}^\alpha(z) \,=\, m_\alpha(x,y,z) z^{2\alpha +1}\]
with the kernel
\[ m_\alpha(x,y,z) \,=\,
\frac{2^{1-2\alpha}\Gamma(\alpha +1)}{\sqrt{\pi}\,\Gamma(\alpha+{1}{2})}\cdot \frac{[(z^2-(x-y)^2)((x+y)^2-z^2)]^{\alpha-1/2}}{(xyz)^{2\alpha}}\cdot 1_{[|x-y|,x+y]}(z).\]
The product formula \eqref{(2.106)} induces a convolution of point 
measures on $\b R_+$ according to  
\[\delta_x\circ_\alpha\delta_y:=\nu_{x,y}^\alpha\,,\]
which in turn allows a  unique bilinear and weakly continuous extension to a probability preserving 
convolution on $M_b(\b R_+).$ It induces the structure of a commutative 
hypergroup on $\b R_+$, called the 
Bessel-Kingman hypergroup of index $\alpha$, c.f \cite{BH}. 
The Dunkl kernel $E_k$ itself  satisfies a product formula of type 
\eqref{(1.1g)}, see \cite{R0} and also \cite{Ros}. Here the 
representing measures $\mu_{x,y}^k$ are not 
positive, but signed 
Borel measures on $\b R$ with uniformly bounded total variation norm. 
For $x,y\not=0$, one has
\[ d\mu_{x,y}^k(z) \,=\, m_{k-1/2}(|x|,|y|,|z|)|z|^{2k}\cdot\frac{1-\sigma_{x,y,z} 
+ \sigma_{z,x,y} + \sigma_{z,y,x}}{2},\]
with
\[ \sigma_{z,x,t} = 
\begin{cases}\displaystyle{\frac{z^2 + x^2 - t^2}{2zx}}
   &\text{if $z,x \not= 0$},\\
     0 \,& \text{else.}
\end{cases}\]
There further exists a unique bilinear and separately weakly continuous convolution $*_k$ on $M_b(\b R)$ such that the convolution of point measures satisfies
$\delta_x*_k\delta_y = \mu_{x,y}^k$. This convolution makes 
$M_b(\b R)$ into 
a commutative Banach-$*$-algebra with unit $\delta_0$, involution 
$\mu^*(A)=\overline{\mu(-A)}$ and  norm $\|\mu\|^\prime = \|L_\mu\|,$ 
the operator $L_\mu$ on $M_b(\b R)$ being defined by 
$L_\mu(\nu) = \mu*_k\nu$.
 The Gelfand transform on $M_b(\b R)$ coincides with the 
corresponding 
Dunkl transform. For details, see \cite{R0}. 
Finally,  the convolution $*_k$ matches the Dunkl-type
generalized translation as defined above. In fact, Lemma \ref{L:Translation}
shows that
$\,\int_{\b R} f\,d(\delta_x*_k\delta_y) = f(x*_k y)$ for all $f\in \mathscr S(\b R)$ and hence, by  a simple density argument, also for all 
$f\in C^\infty(\b R)$.

\section{Positivity of the spherical mean operator}

Following \cite{MT}, we define the Dunkl-type spherical mean operator $f\mapsto M_f$ on $C^\infty(\b R^N)$  by
\[ M_f(x,t):= \frac{1}{d_k} \int_{S^{N-1}} f(x*_k ty) w_k(y) d\sigma(y), 
 \quad (x\in \b R^N, \, t\geq 0).\]
In case $N=1,$ this reduces to
\[ M_f(x,t) \,=\, \frac{1}{2}\bigl(f(x*_k t) + f(x*_k -t)\bigr).\]
Lemma \ref{L:Translation} easily implies 
that $M_f\in C^\infty(\b R^N\times \b R_+)$.
The following is the key result of this paper:

\begin{theorem}\label{P:Pos} The spherical mean operator $f\mapsto M_f$ is positivity-preserving 
on $C^\infty(\b R^N)$, i.e. 
\[ f\geq 0 \,\text{ on }\, \b R^N \quad\text{implies that}\quad M_f\geq 0 \,
\text{ on }\, \b R^N\times \b R_+\,.\]
\end{theorem}

\begin{remark}
The assertion is obvious in the rank-one case: In fact, $M_f$ can then be 
calculated  explicitly in terms of the data given in Section \ref{rank-one}.
We have
\[ M_f(x,t) \,=\, \int_{\b R} f(z) d\sigma_{x,t}^k(z) \]
with 
\[ \sigma_{x,t}^k\,=\, \frac{1}{2}(\mu_{x,t}^k + \mu_{x,-t}^k),\]
which is easily checked to be a probability measure.
\end{remark}

The proof of Theorem \ref{P:Pos} in the higher rank case  is achieved by 
a  
reduction to initial data of the form 
$f(x) = \Gamma_k(s,x,y),$ where $\Gamma_k$   is the Dunkl 
type heat kernel 
\begin{equation}\label{(3.102)} 
\Gamma_k(s,x,y)  = 
\frac{1}{(2s)^{\gamma +N/2}c_k}e^{-(|x|^2 + |y|^2)/4s} 
E_k\bigl(\frac{x}{\sqrt{2s}}\,,
\frac{y}{\sqrt{2s}}\bigr)\end{equation}
$(x,y\in \b R^N, \,s>0),$ 
see \cite{R1}. (Notice that the constant $c_k$ was there differently defined).
We recall from \cite{R1} that  $\Gamma_k$ is strictly positive with
\[ \int_{\b R^N} \Gamma_k(s,x,y) w_k(y) dy = 1 .\]
Moreover, 
\begin{equation}\label{(3.101)}
 \Gamma_k(s,x,y) = \, \frac{1}{c_k^2}\int_{\b R^N} e^{-s|\xi|^2} E_k(-ix,\xi) 
E_k(iy,\xi) w_k(\xi)d\xi\,.
\end{equation}

\begin{lemma} \label{L:heat}  
\begin{enumerate}
\item[\rm{(1)}] For all $x,y,z\in \b R^N$ and $s>0,$ 
\[ \Gamma_k(s,x*_ky,z) \,=\, \frac{1}{c_k^2} \int_{\b R^N} e^{-s|\xi|^2} 
E_k(-iz,\xi)E_k(ix,\xi) E_k(iy,\xi) w_k(\xi) d\xi.\]
In particular, $\Gamma_k(s,x*_ky,z)$ belongs to $\scr S(\b R^N)$ as a
function of each of the arguments $x,y,z$.  
\item[\rm{(2)}] $\displaystyle \int_{\b R^N} \Gamma_k(s,x*_ky,z)   
        w_k(z)dz \,=\, 1.$
\item[\rm{(3)}] If $f\in \mathscr S(\b R^N)$ and $(x,t) \in \b R^N\times \b R_+,$ then 
\[ M_f(x,t) = \lim_{s\downarrow 0} \int_{\b R^N} 
M_{\Gamma_k(s,\,.\,,z)}(x,t)\, f(z) w_k(z)dz.\]
\end{enumerate}
\end{lemma}

\begin{proof} (1) 
By the inversion theorem for the Dunkl transform, formula \eqref{(3.101)}
is equivalent to 
\[ \Gamma_k(s,\,.\,,z)^{\wedge k}(\xi) = c_k^{-1} e^{-s|\xi|^2} 
   E_k(-iz,\xi).\]
Lemma \ref{L:Translation}(3) thus
 implies the stated identity. The rest follows form the
invariance of $\scr S(\b R^N)$ under the Dunkl transform.

(2) In view of  (1), we write
\[ \Gamma_k(s,x*_k y,z) = \widehat{g_{x,y}}^{\,k}(z),
\quad\text{ with }\,\, 
 g_{x,y}(\xi):=  \frac{1}{c_k} e^{-s|\xi|^2} E_k(ix,\xi) E_k (iy,\xi).\]
This gives
\[\int_{\b R^N} \Gamma_k(s,x*_ky,z) w_k(z)dz \,=\, \int_{\b R^N} \widehat 
{g_{x,y}}^{\,k}(z) w_k(z) dz \,=\, c_k g_{x,y}(0) \,=\, 1.\]

(3) By (1) and the definition of the Dunkl transform, one obtains
\begin{align} \int_{\b R^N} &
M_{\Gamma_k(s,\,.\,,z)}(x,t)f(z) w_k(z)dz\notag\\ =&\,\frac{1}{d_k} 
\int_{\b R^N}\Bigl(\int_{S^{N-1}} 
\Gamma_k(s,x*_k ty,z)(x,r) w_k(y) d\sigma(y)\Bigr) 
     f(z)w_k(z) dz\notag\\
=&\, \frac{1}{c_k d_k} \int_{\b R^N}\int_{S^{N-1}} e^{-s|\xi|^2} 
\widehat f^{\,k}(\xi) E_k(ix,\xi)E_k(ity,\xi) w_k(y) d\sigma(y) 
w_k(\xi)d\xi.\notag
\end{align}
Thus, by the dominated convergence theorem,
\begin{align} \lim_{s\downarrow 0} \int_{\b R^N}&
M_{\Gamma_k(s,\,.\,,z)}(x,t) f(z) w_k(z)dz\notag\\ =&\,\, 
      \frac{1}{c_kd_k} \int_{\b R^N}
\int_{S^{N-1}} \widehat f^{\,k}(\xi) E_k(ix,\xi)E_k(ity,\xi)
w_k(y)d\sigma(y)w_k(\xi)d\xi\notag\\
=&\,\, \frac{1}{d_k} \int_{S^{N-1}}f(x*_k ty)w_k(y)d\sigma(y)\,=\,\,
M_f(x,t).\notag \end{align}
\end{proof}

\begin{remark}
A similar calculation shows that for $f\in \scr S(\b R^N)$ 
and  $x,y\in \b R^N$, 
\[ f(x*_k y) \,=\, \lim_{s\downarrow 0} 
\int_{\b R^N} f(z) \Gamma_k(s, x*_k y,z)w_k(z)dz.\]
\end{remark}

\begin{proof}[Proof of Theorem \ref{P:Pos}] We may 
assume that $N\geq 2$  and $\gamma >0,$ hence also $\lambda >0.$ Moreover, it suffices 
to prove the result for $f\in \mathscr S(\b R^N),$ because $\mathscr S(\b R^N)$ is dense in $C^\infty(\b R^N)$ and $M_f(x,t)$ depends continuously on $f$ as a consequence of 
Lemma \ref{L:Translation}(2).
Thus, by part (3) of Lemma \ref{L:heat}, it remains to show that
\begin{equation}\label{(1.9)}
 M_{\Gamma_k(s,\,.\,,z)}(x,t)\, > 0 \quad\text{for all }\, 
s\geq 0,\,z\in \b R^N \text{ and all }\, t\geq 0, \,x\in \b R^N.
\end{equation}
Let us write $M(x,t):= M_{\Gamma_k(s,\,.\,,z)}(x,t)$ for brevity.
Invoking part (1) of Lemma \ref{L:heat} and Corollary \ref{Funk}
leads to 
\begin{align}\label{(1.10)}
 M(x,t) =\,& \frac{1}{d_k}
\int_{S^{N-1}} \Gamma_k(s,x*_kty,z) w_k(y) d\sigma(y)\notag\\
=\, & \frac{1}{c_k^2} \int_{\b R^N} e^{-s|\xi|^2} E_k(-iz,\xi) E_k(ix,\xi)\, j_\lambda (t|\xi|)\, w_k(\xi) d\xi \notag\\
= \, &  \frac{d_k}{c_k^2} \int_0^\infty I(x,z,r)\, e^{-sr^2} j_\lambda (tr) r^{2\lambda +1} dr,
\end{align}
with
\begin{equation}\label{(1.100)}
 I(x,z,r)\,= \,\frac{1}{d_k}
\int_{S^{N-1}} E_k(ix,r\xi) E_k(-iz,r\xi) w_k(\xi) d\sigma(\xi).
\end{equation}
This integral has to be brought in a form from which the positivity of $M(x,t)$ can
be read off. For this, 
we insert the series expansion \eqref{(1.2)} for $E_k(ix,r\xi)$, 
thus obtaining
\begin{align}\label{(1.11)} I(x,z,r)\, =&\,\, \frac{1}{d_k}\sum_{n=0}^\infty 
\frac{\Gamma(\lambda +1)}{2^n\,\Gamma(n+\lambda+1)} j_{n+\lambda} (r|x|)\cdot
\notag\\
&\qquad\qquad\cdot\int_{S^{N-1}} P_n^k(irx,\xi) E_k(-irz,\xi) w_k(\xi) d\sigma(\xi)\notag\\ 
=&\,\, \sum_{n=0}^\infty  \Bigl(
 \frac{\Gamma(\lambda +1)}{2^n\,\Gamma(n+\lambda+1)}\Bigr)^2 
 j_{n+\lambda}(r|x|) j_{n+\lambda}(r|z|) P_n^k(irx, -irz),
\end{align}
where for the second identity, again Corollary \ref{Funk} was used.
In case $x,z\not=0$, the homogeneity of $P_n^k$ allows one to write
\begin{equation}
\label{(1.12)}
P_n^k(irx, -irz) \,=\, \frac{(n+\lambda)(2\lambda)_n}{\lambda\cdot n!}(r^2|x||z|)^n \,
        V_k\widetilde C_n^\lambda\Bigl
(\Big< \frac{x}{|x|},\,.\,\Big>\Bigr)\Bigl(\frac{z}{|z|}\Bigr).
\end{equation}
We now employ a well-known degenerate version of the 
 addition theorem for Gegenbauer polynomials (see 
\cite{A}, (4.36)): For all $s,t,\theta\in \b R$,
\begin{align}
j_\lambda\bigl(\sqrt{s^2 + t^2 - 2st \cos\theta}\,&\bigr)\notag\\
=\,
\sum_{n=0}^\infty \frac{(n+\lambda)(2\lambda)_n}{\lambda\cdot n!} \Bigl(&\frac{\Gamma(\lambda +1)}
{2^n\,\Gamma(n+\lambda +1)}\Bigr)^2 (st)^n j_{n+\lambda}(s)
j_{n+\lambda}(t)\, \widetilde C_n^\lambda(\cos\theta).\notag
\end{align}
This series converges uniformly with respect to $\theta\in \b R$. 
Combining \eqref{(1.11)} and \eqref{(1.12)} and recalling the Laplace representation \eqref{(1.40a)} 
for $V_k$, we now see that for all $x,z\not=0,$
\begin{align}\label{(1.101)}
I(x,z,r)\,=\,& \int_{\b R^N} j_\lambda \,\bigl(r\sqrt{|x|^2 + |z|^2 
-2 |z| \langle x,\eta\rangle}\,\bigr)\, d\mu_{z/|z|}^k(\eta)\notag\\ 
=\,& \int_{\b R^N} j_\lambda \,\bigl(r\sqrt{|x|^2 + |z|^2 
-2  \langle x,\eta\rangle}\,\bigr)\, d\mu_{z}^k(\eta);
\end{align}
The second identity follows from the dilation 
 equivariance of $\mu_z^k$:  $\mu_{rz}^k(A) = \mu_z^k(r^{-1}A)$ for all $r>0, \, A\in \mathcal B(\b R^N),$ c.f. \cite{R2}. 
By virtue of Corollary \ref{Funk}, \eqref{(1.101)} remains true if $x=0$ or $z=0$.
  For abbreviation, put
\[ v_z(\eta):= \sqrt{|x|^2 + |z|^2 -2\langle x,\eta\rangle}.\]
Then by \eqref{(1.10)}, the product formula \eqref{(2.106)}
 for the Bessel functions $j_\lambda$ and eq. 11.4.29 of \cite{AS},
 we arrive at
\begin{align}
 M(x,t) \,=& \, \frac{d_k}{c_k^2} \int_{\b R^N}\int_0^\infty j_\lambda(rv_z(\eta))\,
j_\lambda(rt)\, e^{-sr^2}
r^{2\lambda+1}dr\, d\mu_z^k(\eta)\notag\\
=&\, \frac{d_k}{c_k^2}\int_{\b R^N} \int_0^\infty \Bigl(\int_0^\infty j_\lambda(ru) e^{-sr^2}
r^{2\lambda+1} dr\Bigr)\,
d\nu^\lambda_{v_z(\eta),t}(u)\,d\mu_z^k(\eta)\notag\\
=&\,\frac{d_k}{c_k^2}\frac{\Gamma(\lambda+1)}{2s^{\lambda+1}} \int_{\b R^N} \Bigl(\int_0^\infty
  e^{-u^2/4s}d\nu_{v_z(\eta),t}^\lambda(u)\Bigr)d\mu_z^k(\eta),
\notag
\end{align}
which is obviously non-negative. This 
finishes the proof.
\end{proof}

\section{A positive radial product formula for the Dunkl kernel}

\subsection{Statement of the main result}

\noindent
Consider  $f(x) = E_k(ix,z)$ with 
$z\in \b R^N$. Then
\begin{equation}\label{(1.210)}
 M_f(x,t) = E_k(ix,z) j_\lambda(t|z|).
\end{equation}
Indeed, we have
\[ \frac{1}{d_k}\int_{S^{N-1}} E_k(ix,z) w_k(z)d\sigma(z)\,=\, j_\lambda(|x|),\]
which follows from Corollary \ref{Funk} if $N\geq 2$ and is  obviously also true if $N=1$.
In view of  \eqref{(1.5b)}, this gives \eqref{(1.210)}. 
Theorem \ref{P:Pos} easily implies the existence of the representing measures stated 
in the next theorem, which is the main result of this paper:

\begin{theorem}\label{T:Main} For each $x\in \b R^N$ and $t\in \b R_+$ there exists a unique compactly supported probability measure
$\sigma^k_{x,t} \in M^1(\b R^N)$ such that
\begin{equation}\label{(1.9a)}
 E_k(ix,y)\,j_\lambda(t|y|) \,=\, 
\int_{\b R^N} E_k(i\xi, y) d\sigma^k_{x,t}(\xi) 
   \quad \text{for all } \,y\in \b R^N.
\end{equation}
It represents the spherical mean operator $f\mapsto M_f$ in the sense that
\begin{equation}\label{(1.9b)} 
 M_f(x,t) \,=\, \int_{\b R^N} f\,d\sigma_{x,t}^k \quad\text{for all } \,f\in 
  C^\infty(\b R^N).
\end{equation}
The measure $\sigma^k_{x,t}$ satisfies
\[\text{supp}\,\sigma^k_{x,t}
\subseteq \,
\bigcup_{g\in G}\{\xi\in \b R^N: |\xi-gx|\leq t\},\]
and the mapping $(x,t) \mapsto \sigma_{x,t}^k$ is  continuous with respect to the weak topology on $M^1(\b R^N).$
 Moreover,  
\[\sigma^k_{gx,t}(A) = \sigma_{x,t}^k(g^{-1}(A))\quad\text{and }\,
\sigma_{rx,rt}^k(A) = \sigma_{x,t}^k
(r^{-1}A)\]
for all $g\in G, r>0,$ and all Borel sets $A\in \mathcal B(\b R^N)$.  
\end{theorem}

\begin{remarks}
\begin{enumerate}
\item[\rm{1.}]
The above result also gives a natural extension of the spherical mean operator
at hand: namely for  $f\in C_b(\b R^N),$  we may define $M_f\in C_b(\b R^N\times
\b R_+)$ by
\[ M_f(x,t):= \int_{\b R^N} fd\sigma_{x,t}^k\,.\]
\item[\rm{2.}]
In the rank-one case, the support of $\sigma_{x,t}^k$ with $t\not=0$ is given by
\[ \text{supp}\,\sigma_{x,t}^k\,=\, 
\big[-|x|-t, -||x|-t|\big]\cup \big[||x|-t|, |x|+t\big].\]
This 
illustrates that the complete $G$-orbit of $x$ has to be taken into account in the 
description of $\text{supp}\,\sigma_{x,t}^k.$
\item[\rm{3.}] From part (2) of Lemma \ref{L:Translation} a weaker statement  
on the support of $\sigma_{x,t}^k$ is immediate, namely $\,\text{supp}\, 
\sigma_{x,t}^k \subseteq \{\xi\in \b R^N: |\xi|\leq |x|+t\}.$
\end{enumerate}
\end{remarks}

In order to derive the stated properties of the support of 
$\sigma_{x,t}^k$ we use an approach via a Darboux type initial value problem 
satisfied
by the spherical mean $M_f$.  For $\alpha\geq -1/2,$ denote by $A_\alpha$ the 
singular Sturm-Liouville operator 
\[ A_\alpha\,:=\, \partial_t^2 + \frac{2\alpha+1}{t}\partial_t\]
on $\b R_+$. The Bessel functions $t\mapsto j_\alpha(tz),\, z\in \b C,$ are (up to normalization) the unique even and
analytic eigenfunctions of this operator, satisfying
\[ A_\alpha^t j_\alpha(tz) \,=\, -z^2 j_\alpha(tz);\]
here the uppercase index indicates the relevant variable.

For 
$f\in \mathscr S(\b R^N),$ we may write
\begin{equation}\label{(1.21)}
 M_f(x,t) = \,\frac{1}{c_k}\int_{\b R^N} \widehat f^{\,k}(\xi) E_k(x,i\xi) j_\lambda(t|\xi|)
    w_k(\xi)d\xi.
\end{equation}
In this case, a direct calculation verifies that $u=M_f$ solves the 
following 
initial value problem for the Darboux-type differential-reflection operator 
$\Delta_k^x - A_\lambda^t$:
\begin{align} \label{(1.8)}
 &(\Delta_k^x - A_\lambda^t) u \,=\, 0 \quad\text{in }\, 
 \b R^N\times \b R_+\,;\\ 
 &u(x,0) = f(x), \,\, u_t (x,0) = 0 \quad\text{for all } x\in \b R^N.\notag
\end{align}
(To obtain the initial data, one has to 
use the inversion theorem for the Dunkl transform as well as the fact 
that $M_f$ is even 
with respect to $t$). 
In \cite{MT}, Prop. 5.2 it is shown that for arbitrary $f\in C^\infty(\b R^N),$
$M_f$ is in fact the unique $C^\infty$-solution of \eqref{(1.8)}.
Instead of studying the domain of dependence for  the above Darboux-type 
operator directly, we shall consider  the corresponding reflection invariant wave operator, 
which is easier to handle:

\subsection{Domain of dependence for wave operators 
related with reflection
      groups}

This section is devoted to the study of the Dunkl-type wave operator $\,L_k -\partial_t^2$, where  
$L_k$ denotes the reflection invariant part of the Dunkl Laplacian
$\Delta_k$, c.f.  \eqref{(2.108)}.
Notice that in the rank-one case,  $L_k$ coincides with the Sturm-Liouville 
operator $A_{k-1/2}$.  
It will be important for the following that
$L_k$ can  be written in divergence form, 
\begin{equation}\label{(1.3a)}
 L_k \,=\, \frac{1}{w_k(x)} \sum_{i=1}^N 
\partial_{x_i}\bigl(w_k(x)\partial_{x_i}\bigr).\end{equation}
Let us now turn to the Dunkl-type
wave equation on $\b R^N\times\b R_+$ 
associated with $G$ and $k$,
\begin{equation}\label{(1.4)} 
 (L_k-\partial_t^2) u\,=\,0.
\end{equation} 
The following uniqueness result for solutions of this 
equation is in close analogy
to well-known  facts in the classical case $k=0$.

\begin{lemma}
Suppose that $u$ is a real-valued $C^2$-solution of \eqref{(1.4)},
given in  the truncated cone
 \[C(x_0,t_0)= 
\{(x,t)\in \b R^N\times \b R_+\,: |x-x_0|\leq t_0-t\}.\]
Define the energy
of $u$ at time $t$ within this cone by
\[ E(t) := \frac{1}{2}\int_{B_t}  \bigl(u_t^2 + |\nabla_{\!x}u(x,t)|^2\bigr) 
w_k(x) dx,\]
 where the integration is taken over the  level set 
\[B_t =\{x\in \b R^N: |x-x_0|\leq t_0 -t\}.\] 
Then $\,E^\prime \leq 0$ on $\b R_+$. 
\end{lemma}

\begin{proof} The chain rule gives
\[ E^\prime(t) \,=\, I(t) -\frac{1}{2} \int_{\partial B_t} \bigl(u_t(x,t)^2 + 
|\nabla_x u(x,t)|^2\bigr) w_k(x) dx\]
with
\begin{align}
 I(t) =\, \int_{B_t}& \bigl( u_t(x,t) u_{tt}(x,t) + \sum_{j} u_{x_j}(x,t)u_{t,x_j}
(x,t)\bigr) 
              w_k(x)dx\notag\\
     =\, \int_{B_t}& \bigl( u_{tt}(x,t) - L_k^x u(x,t)\bigr) 
u_t(x,t) w_k(x)dx \notag\\
 &+\,\int_{B_t} \sum_{j=1}^N \bigl( u_t(x,t) u_{x_j}(x,t) w_k(x)\bigr)_{x_j} dx
 \notag\\
  =\,\int_{\partial B_t}& u_t(x,t)\, \partial_\nu u(x,t)\,w_k(x) 
  d\sigma(x),\notag
\end{align}
where the second identity follows from \eqref{(1.3a)} and the last one
 from the divergence theorem ($\nu$ is the exterior unit normal to $B_t$ 
and $d\sigma$ the Lebesgue surface measure on $\partial B_t$). By the 
Cauchy-Schwarz
 inequality,
\[ u_t \partial_\nu u\,\leq \,\frac{1}{2} \bigl(u_t^2 + 
|\nabla_xu|^2\bigr).\]
This implies the assertion.
\end{proof}

\noindent
The following uniqueness result is an immediate consequence:

\begin{theorem}\label{T:depend}
Suppose that $u$ is a $C^2$-solution of the wave equation $(L_k- \partial_t^2)u=0$, 
 defined in the cone 
$C(x_0,t_0)$ and satisfying 
\[  u_t(x,0) = u(x,0) = 0 \quad\text{for all } x\in \b R^N
\text{ with } |x-x_0| \leq t_0.\]
 Then $u$ vanishes in  $C(x_0,t_0)$. 
\end{theorem} 

\begin{proof} We may assume that $u$ is real-valued. As $E(t) \geq 0$ and $E(0) =0$,
 the lemma shows that
$E(t) =0$ for  $0\leq t\leq t_0$, and hence $u_t =0$ and $\nabla_{\! x}u = 0$ in 
$C(x_0,t_0)$. This implies the assertion.
\end{proof}

We mention that there exists a thorough study of wave operators (and more general hyperbolic operators) 
related with root systems especially in the case of integer-valued  multiplicity functions, see \cite{B} and the references cited there.  In particular,  for integer-valued  $k$ and 
odd dimensions $N$ satisfying $N\geq 2\gamma+3$, 
the wave operator $\,L_k-\partial_t^2\,$ in fact satisfies 
the strong Huygens principle. This means 
that a solution at $(x_0,t_0)$ depends only on the data in an infinitesimal
neighborhood of the surface of the propagation cone with 
vertex at $(x_0,t_0)$.

\subsection{Proof of the main result}\label{S:Pos}

For the proof of the properties of the representing measures in Theorem \ref{T:Main}, we have  to relate the Darboux-type equation
\eqref{(1.8)} with the hyperbolic equation 
\eqref{(1.4)}. For this, we   involve the Riemann-Liouville
transform with parameter $\alpha >-1/2$ on $\b R_+$. It is given by
\begin{equation}\label{(1.8a)}
 \mathcal R_\alpha f(t)\,=\, 
\frac{2\Gamma(\alpha+1)}{\Gamma(1/2)\Gamma(\alpha +1/2)} 
\int_0^1 f(st) (1-s^2)^{\alpha-1/2}ds 
\end{equation}
for $f\in C^\infty(\b R_+),$ see \cite{T}.
The operator $\mathcal R_\alpha$  satisfies the intertwining property
\[ A_\alpha \mathcal R_\alpha\,=\,\mathcal R_\alpha\, \frac{d^2}{dt^2}.\]  
Notice that  the rank-one intertwining operator
$V_k$, when restricted to even functions, just coincides with 
$\mathcal R_{k-1/2}$. 
According to \cite{T}, $\mathcal R_\alpha$ 
is a topological isomorphism of  $C^\infty(\b R_+)$ with respect to
the usual Fr\'echet-topology,
which is induced from the identification of  $C^\infty(\b R_+)$ with 
$\{f\in C^\infty(\b R): f(-t) = f(t)\}$.

\begin{proof}[Proof of Theorem \ref{T:Main}]
For fixed $x\in \b R^N$ and $t\geq 0$, consider the linear functional
\[\Phi_{x,t}: \,f\mapsto M_f(x,t),\]
which is positivity-preserving on $C^\infty(\b R^N)$ according to 
Theorem \ref{P:Pos}.
Moreover, $\Phi_{x,t}(1) = 1$.
It follows that $\Phi_{x,t}$ is  represented by 
a compactly supported probability  measure $\sigma_{x,t}^k \in M^1(\b R^N)$ (c.f. Theorem 2.1.7 of \cite{Ho}). This implies  statement \eqref{(1.9b)}  and also \eqref{(1.9a)} in view of relation \eqref{(1.210)}. 
Notice that
\begin{equation}
\label{(2.10)} 
E_k(ix,z)j_\lambda(t|z|) \,=\, (\sigma_{x,t}^k)^{\wedge k}(-z).
\end{equation}
Thus 
the uniqueness of $\sigma_{x,t}^k$ follows from the injectivity of 
the Dunkl  transform on $M^1(\b R^N)$, c.f. Propos. \ref{Dunkltrafo}.
In order to check the weak continuity of $(x,t) \mapsto \sigma_{x,t}^k$, 
take a sequence $(x_n, t_n)_{n\in \b N}\subset \b R^N\times \b R_+$ with $\lim_{n\to\infty}(x_n,t_n) = (x_0,t_0)$. Then \eqref{(2.10)} implies that $(\sigma_{x_n,t_n}^k)^{\wedge k}\, \to \,(\sigma_{x_n,t_n}^k)^{\wedge k}$ pointwise on $\b R^N$. 
L\'evy's continuity theorem for the Dunkl transform (Lemma \ref{Dunkltrafo} (5))
now yields that the  
 $\sigma_{x_n,t_n}^k$ converge weakly to $\sigma_{x_0,t_0}^k.$
Further, the 
claimed  transformation  properties of $\sigma^k_{x,t}$  
are immediate consequences of the invariance properties \eqref{(2.102)}
 of the kernel $E_k$.
It remains to analyse the support  of $\sigma_{x,t}^k$. For this, 
we shall employ Theorem \ref{T:depend}. We therefore first  
reduce the general situation
to the group-invariant case. We introduce the group means
\[ \widetilde\sigma_{x,t}^{\,k}\,:=\, \frac{1}{|G|}\sum_{g\in G}  
\sigma_{gx,t}^k\,=\, \frac{1}{|G|}\sum_{g\in G}  
\sigma_{x,t}^k\circ g^{-1}.  \]
As $\sigma_{x,t}^k$ is positive, it is enough to show that 
$\widetilde\sigma_{x,t}^{\,k}$
is supported in 
\[ K(x,t):= \bigcup_{g\in G}\{\xi\in \b R^N: |\xi-gx|\leq t\},\]
as stated for $\sigma_{x,t}^k$. 
For this in turn, it suffices to show that
\[ \int_{\b R^N} f\,d\widetilde\sigma^k_{x,t}\,= 0 \]
for all $G$-invariant $f\in\mathscr S(\b R^N)$ whose support does not
intersect $K(x,t)$. So suppose  that 
$f\in \mathscr S(\b R^N)$ is $G$-invariant.
Then according to Prop. 2.4. of \cite{RV}, $\widehat f^{\,k}$ is also $G$-invariant, and \eqref{(1.21)} shows 
that $x\mapsto M_f(x,t)$ is $G$-invariant as well. Put
 $ u_f(x,t):= 
(\mathcal R_\lambda^t)^{-1}M_f(x,t),$ which is still $G$-invariant with respect to $x$. According to \eqref{(1.8)} and 
the intertwining property of the Riemann-Liouville transform, $ u=
  u_f$ belongs to 
$C^\infty(\b R^N\times\b R_+)$
 and solves the initial value problem
\begin{align}\label{(1.40)} 
 &(L_k^x  -  \partial_t^2) u \,=\, 0 \quad\text{in }\, \b R^N\times 
\b R_+;\notag\\ 
 &u(x,0) = f(x), \,\, u_t (x,0) = 0 \quad\text{for all } x\in \b R^N.
\end{align}
Now suppose in addition that $\,\text{supp}\,f\cap K(x,t) = \emptyset$. 
Then Theorem \ref{T:depend} implies that $\, u_f(x,s) =0$ for all $0\leq s\leq t$.
From  the explicit form \eqref{(1.8a)} of the Riemann-Liouville transform 
$\mathcal R_\lambda$ we further deduce that  
\[ \int_{\b R^N} f\,d\widetilde\sigma_{x,t}^k\,=\, M_f(x,t) \,=\, 
\mathcal R_\lambda^t u_f(x,t)\,=\,0\]
as claimed. This completes the proof of Theorem \ref{T:Main}.
\end{proof}

\section{Positive translation of radial functions}

In this section, we derive a slightly weaker variant of Theorem \ref{T:Main}.
For its formulation, we introduce some additional notation:
If $\mathcal F$ is a space of $\b C$-valued functions on $\b R^N,$ 
 denote by 
\[\mathcal F_{rad}:= \{f\in \mathcal F: f\circ A = f 
\quad\text{for all }A\in O(N,\b R)\}\]
 the subspace of those $f\in \mathcal F$ which are radial. For $f\in \mathcal F_{rad}$ there exists a unique function $\widetilde f: \b R_+ \to \b C$ such that
$f(x) = \widetilde f(|x|)$ for all $x\in \b R^N.$ 
Similar, if $\mathcal M$ is a space of Borel measures on $\b R^N,$ then 
\[ \mathcal M_{rad}:= \{\mu\in \mathcal M:  \mu\circ A = \mu \quad \text{for all } 
A\in O(N,\b R)\}\] 
denotes the subspace of radial measures from $\mathcal M$. 


\begin{theorem}\label{rad}
For each $x,y\in \b R^N$  there exists a unique compactly supported, 
radial probability 
measure $\rho_{x,y}^k\in M^1_{rad}(\b R^N)$ such that for all 
$f\in C^\infty_{rad}(\b R^N),$
\begin{equation}\label{(1.30)}
 f(x*_ky) = \int_{\b R^N} fd\rho_{x,y}^k\,.
\end{equation}
The support of $\rho_{x,y}^k$ is contained in 
\[\{\xi\in \b R^N: \,
\min_{g\in G}|x+gy| \leq |\xi| \leq \max_{g\in G}|x+gy|\}.\]
 In particular, if $\,0\in \text{supp}\,
\rho_{x,y}^k,$ then  the $G$-orbits of $x$ and $-y$ coincide.
\end{theorem}

\begin{proof}
Take $f\in \mathscr S_{rad}(\b R^N)$. Then according to Prop 2.4. of \cite{RV}, 
$\widehat f^{\,k}\in \mathscr S_{rad}(\b R^N)$ with $\,\widehat f^{\,k}(\xi) =
\mathcal H^{\lambda}(\widetilde f\,)(|\xi|).$ Here $\mathcal H^{\lambda}$
stands for  the Hankel transform of index $\lambda$ on $L^1(\b R_+\,, r^{2\lambda+1}dr)$,
given  by
\[ \mathcal H^\lambda(g)(s) = \frac{1}{2^\lambda\Gamma(\lambda +1)}\int_0^\infty
 g(r)j_\lambda(rs)r^{2\lambda+1} dr.\]
Employing  Lemma \ref{L:Translation} (3), relations  \eqref{(1.100)}, 
\eqref{(1.101)} as well as the 
inversion theorem for the Hankel transform, 
 one obtains
\begin{align}
 f(x*_ky)  \,=&\, \frac{1}{c_k} \int_0^\infty 
  \mathcal H^\lambda (\widetilde f\,)(r)\int_{S^{N-1}} E_k(ix,r\xi) E_k(iy, r\xi)w_k(\xi) d\sigma(\xi) \,r^{2\lambda+1} dr \notag\\
=\,& \frac{d_k}{c_k} \int_0^\infty \mathcal H^\lambda (\widetilde f)(r) 
I(-x,y,r) r^{2\lambda+1}dr \notag\\
=\,& \frac{1}{2^\lambda\Gamma(\lambda+1)} \int_{\b R^N}
\int_0^\infty 
          \mathcal H^\lambda (\widetilde f\,)(r) j_\lambda(r\sqrt{|x|^2 + |y|^2 +
 2\langle x,\eta\rangle})r^{2\lambda +1} dr\, d\mu_y^k(\eta) \notag\\
=\,& \int_{C(y)} \widetilde f\bigl(\sqrt{|x|^2 + |y|^2 + 
           2\langle x,\eta\rangle}\bigr) d\mu_y^k(\eta);\notag
\end{align}
here for  the second identity, it was  used that $I(x,-y,r) = I(-x,y,r)$. 
This shows that  $f(x*_ky) \geq 0$ if $f\geq 0$ on $\b R^N$. 
As $\mu_y^k$ is a compactly supported probability measure, it also follows that
\[ \sup\{|f(x*_ky)|: f\in \mathscr S(\b R^N), \|f\|_\infty\leq 1\} \,=\, 1.\]
Similar to the proof of Theorem \ref{T:Main}, we proceed by considering the 
linear functional 
$\, \Psi_{x,y}\!:   f\mapsto f(x*_ky),\,$
which is positive and bounded with norm $\|\Psi_{x,y}\| = 1$
on the dense subspace $\mathscr S_{rad}(\b R^N)$ of 
$(C_{0,rad}(\b R^N), \|\,.\,\|_\infty)$. Its continuous extension to $C_{0,rad}(\b R^N)$ is 
 a probability measure $\rho_{x,y}^k \in M_{rad}^1(\b R^N)$ which satisfies  \eqref{(1.30)}.
Finally, the statement concerning the support of $\rho_{x,y}^k$ is immediate from the implication
\[ \eta\in C(y) \,\Longrightarrow \,\,\min_{g\in G} |x+gy| \leq
   \sqrt{|x|^2 + |y|^2 + 2\langle x,\eta\rangle} \,\leq\, \max_{g\in G} |x+gy|.\]
\end{proof}

The information on the support of the measures $\rho_{x,y}^k$ can in turn be used to
obtain a sharper statement on the supports of the representing measures for the
spherical mean operator in Theorem \ref{T:Main}:

\begin{corollary} The measures $\sigma_{x,t}^k$ from Theorem \ref{T:Main} satisfy
\[ \text{supp}\,\sigma_{x,t}^k \,\subseteq\, 
 \{\xi\in \b R^N: |\xi|\geq 
 ||x|-t|\}.\]
In particular, if $\,0\in \text{supp}\,\sigma_{x,t}^k,$ then $|x|=t.$ 
\end{corollary}

\begin{proof}
Suppose in the contrary that $\text{supp}\, \sigma_{x,t}^k\nsubseteq 
\{\xi\in  \b R^N: |\xi|\geq ||x|-t|\}.\,$
 Then there exists some $f\in 
\mathscr S_{rad}(\b R^N)$ with $f\geq 0,$ 
\begin{equation}\label{(4.4)} 
\text{supp}f \cap 
\{\xi\in \b R^N: |\xi|\geq ||x|-t|\} \,=\,\emptyset,
\end{equation}
and such that $M_f(x,t) >0.$ But then 
 $\eta\mapsto f(x*_kt\eta)$ is not identically zero on $S^{N-1}$. 
 In view of Theorem \ref{rad}, this is  
a contradiction to \eqref{(4.4)}. 
\end{proof}

\section{An application: Semigroups of $k$-invariant Markov kernels with radial distributions}
\label{S:Markov}

There is a concept of  homogeneity  for Markov processes
on $\b R^N,$ called $k$-in\-vari\-ance, which generalizes the classical 
notion of processes with independent,
stationary increments to the Dunkl setting. This was introduced and studied in some 
detail in \cite{RV}. The most important example
is a generalization of Brownian motion, with the transition probabilities 
given in terms of the generalized heat kernel $\Gamma_k$, 
c.f. example \ref{E:ex}(1) below. A Dunkl-type Brownian motion 
is a Feller process and therfore admits a version with c\`adl\`ag-paths. 
This version gives,
 after symmetrization with respect to the 
underlying reflection group, a diffusion on the Weyl chambers. 
$k$-invariant Markov processes are constructed from semigroups of Markov kernels
which are $k$-invariant in the following sense:

\begin{definition}\!(\cite{RV}) 
A Markov kernel $P: \b R^N\times \mathcal B(\b R^N) \to [0,1]$ is called 
$k$-invariant, if
\[ P(x,\,.\,)^{\wedge k}\!(\xi) \,=\, P(0,\,.\,)^{\wedge k}\!(\xi)\cdot
E_k(-ix,\xi) \quad\text{for all } x,\xi \in \b R^N.\]
Here $P(x,\,.\,)$ is regarded as a probability measure on $\b R^N.$ 
\end{definition}

If $k=0$, the $k$-invariant Markov kernels are exactly those which are translation invariant, which means that they 
 satisfy $P(x+x_0, A+x_0) = P(x,A)$ for all $x,x_0\in \b R^N$ 
and all $A\in \mathcal B(\b R^N).$ Equivalently, the translation invariant Markov kernels 
(with $k=0$) are those of the form 
 $P(x,A):= \delta_x*\mu(A)$ with a probability 
measure $\mu\in M^1(\b R^N)$;  here $\ast$ denotes the usual group convolution. 
If $k>0,$ then for given $\mu\in M^1(\b R^N)$
there usually  exists no $k$-invariant 
Markov kernel such that $P(0,\,.\,) = \mu$. This is due to the fact that
the associated generalized translation on $\b R^N$ cannot be expected to be
probability-preserving (and is definitely not in the rank-one case).  

However, 
our results allow us to define associated $k$-invariant Markov kernels 
for all  measures $\mu\in M^1(\b R^N)$ which belong to the class 
\[M_k^1(\b R^N):= \{\mu\in M^1(\b R^N): 
w_k^{-1}(x)d\mu(x) \,\,\text{ is radial.}\} \]
This leads to a considerable variety of 
$k$-invariant Markov processes beyond those discussed in \cite{RV}.
Concerning the definition of $M_k^1(\b R^N),$ one should notice 
 that for  $\mu\in M_k^1(\b R^N),$ the measure $w_k^{-1}d\mu\,$ need not be a 
Radon measure.

\begin{proposition}\label{T:Markov}
Let $\mu\in M_k^1(\b R^N).$ 
Then for each $x\in \b R^N,$ there exists a unique probability measure 
$\delta_x\ast_k\mu\in M^1(\b R^N)$ such that
\begin{equation}\label{(4.2)}
 (\delta_x \ast_k \mu)(f)= \int_{\b R^N} f(x\ast_ky) d\mu(y)  \quad
            \text{for all }\,f\in \mathscr S(\b R^N).
\end{equation}
It satisfies
\begin{equation}\label{(4.2a)}
 (\delta_x\ast_k\mu)(f) \,=\, \int_{\b R^N} M_f(x,|y|) d\mu(y) 
           \quad\text{for all }\,f\in C_b(\b R^N).\end{equation}
\end{proposition}

\begin{proof} Recall first that $y\mapsto f(x\ast_k y)$ belongs to $\mathscr S(\b R^N)$ 
for  $f\in \mathscr S(\b R^N),$
whence the integral on the right side of \eqref{(4.2)} is well-defined.
If $\mu\in M^1_k(\b R^N)$ then for 
$f\in \mathscr S(\b R^N)$, 
\begin{align}
\int_{\b R^N} f(x\ast_ky)d\mu(y)\,= &\,
  \int_{\b R^N}\frac{1}{d_0} \Bigl(\int_{S^{N-1}} f(x*_k|y|\xi) 
    w_k(|y|\xi)d\sigma(\xi)\Bigr) w_k^{-1}(y)d\mu(y)\notag\\
=&\,\frac{d_k}{d_0} \int_{\b R^N} M_f(x,|y|)|y|^{2k} w_k^{-1}(y) d\mu(y)\,
=\, \int_{\b R^N} M_f(x,|y|) d\mu(y).\notag
\end{align}
This shows that the functional $\, f\mapsto 
\int_{\b R^N} f(x\ast_ky)d\mu(y)\,$ extends uniquely to a positive bounded Borel
measure $\delta_x\ast_k\mu\in M_b(\b R^N)$ and that \eqref{(4.2a)}
 is  satisfied.
As $M_f = 1$ for $f= 1$,  
$\delta_x\ast_k\mu$ is in fact a probability measure.
\end{proof}

\begin{theorem}\label{T:Markov2}
\begin{enumerate}
\item[\rm{(i)}]
If $\mu\in M_k^1(\b R^N),$ then $\,P(x,A):= \delta_x*_k\mu(A)$ 
defines a $k$-invariant Markov kernel on $\b R^N$ 
with $P(0,\,.\,) = \mu.$ 
\item[\rm{(ii)}]
Conversely, if $P$ is a  $k$-invariant Markov kernel
with $P(0,\,.\,)\in M_k^1(\b R^N)$, then $P(x,\,.\,) = \delta_x*_k P(0,\,.\,).$
\end{enumerate}
\end{theorem}

\begin{proof}
For (i), notice first that
 $\delta_0\ast_k\mu = \mu.$ Moreover, by the above proposition and  
equation 
\eqref{(1.210)},
\begin{align}
 (\delta_x\ast_k\mu)^{\wedge k}\!(\xi) \,=\,& \int_{\b R^N} E_k(y,-i\xi)d(\delta_x\ast_k\mu)(y)\,=\, 
   \int_{\b R^N} M_{E_k(\,.\,,-i\xi)}(x,|y|)d\mu(y)\notag\\
=&\,E_k(x,-i\xi) \int_{\b R^N} j_\lambda(|y||\xi|) d\mu(y).\notag
\end{align}
Specializing to $x=0,$ this gives 
\begin{equation}\label{(5.100)} 
\widehat \mu^{\,k} (\xi)\,=\, \int_{\b R^N} j_\lambda(|y||\xi|) d\mu(y)
\end{equation}
and thus
\[(\delta_x\ast_k\mu)^{\wedge k}\!(\xi)\,
=\,E_k(-ix,\xi)\cdot\widehat\mu^{\,k}(\xi),\]
which proves that $P$ is $k$-invariant. By L\'evy's continuity theorem for the Dunkl
transfrom, this 
also 
ensures that
$x\mapsto \delta_x*_k\mu$ is continuous with respect
to the weak topology on $M^1(\b R^N).$ This
 easily 
implies that $x\mapsto P(x,A)$ is measurable for any $A\in \mathcal B(\b R^N).$ 
Thus the proof of (i) is complete. Part (ii) follows immediately from the 
definition of $k$-invariance and the injectivity of the Dunkl transform.
\end{proof}

We also introduce a convolution product for measures from $M_k^1(\b R^N)$:

\begin{definition}
Let $\mu,\,\nu\in M_k^1(\b R^N)$. Then $\mu*_k\nu\in M^1(\b R^N)$ is defined by
\[ \mu*_k\nu(f):= \int_{\b R^N} \delta_x*_k\nu(f)d\mu(x)\,=\, 
\int_{\b R^N}\int_{\b R^N} f(x*_ky) d\mu(x)d\nu(y), \quad f\in \mathscr S(\b R^N).
\] 
\end{definition}

\noindent
Notice that 
\begin{equation}\label{(5.101)}
(\mu*_k\nu)^{\wedge k} =\, \widehat\mu^{\, k}\widehat\nu^{\, k} 
           \quad\text{for all }\mu,\,\nu\in M_k^1(\b R^N).
\end{equation}
There is a close relationship between $*_k$ on $M_k^1(\b R^N)$ and the 
 convolution $\circ_\lambda$ of the 
Bessel-Kingman hypergroup of index $\lambda= \gamma+N/2-1$ on $\b R_+$ 
(c.f. Section \ref{rank-one}). We shall in particular  obtain from this connection
 that $M_k^1(\b R^N)$ is
closed w.r.t. $*_k$. Let us start with some additional notation:

We denote by
  $p(\mu)\in M_b(\b R_+)$ the image measure of $\mu\in M_b(\b R^N)$ under the mapping 
$x\mapsto |x|,$ characterized by
\[ \int_{\b R^N} f(|x|)d\mu(x)\,=\, \int_0^\infty f(r) dp(\mu)(r) \quad \text{for all } f\in C_b(\b R_+).\] 
Further, we write $\mathcal H^\lambda$ for the Hankel transform of index 
$\lambda$ on $M_b(\b R_+),$ i.e. 
\[ \mathcal H^\lambda(\sigma) (r) =\, \int_0^\infty j_\lambda(rt)d\sigma(t), \quad \sigma\in M_b(\b R_+).\]
In these terms, \eqref{(5.100)} just states that for $\mu\in  M_k^1(\b R^N),$
\begin{equation}\label{(5.110)}
\widehat \mu^k(\xi)\,=\,\mathcal H^\lambda\bigl(p(\mu)\bigr)(|\xi|).
\end{equation}
We also recall from hypergroup theory that the Hankel transfrom $\mathcal H^\lambda$ is 
multiplicative w.r.t. the convolution $\circ_\lambda$ on $M_b(\b R_+),$ i.e. 
\[ \mathcal H^\lambda(\sigma\circ_\lambda\tau) = \, 
\mathcal H^\lambda(\sigma) \mathcal H^\lambda(\tau)\quad\text{for all }\,
\sigma,\,\tau\in M_b(\b R_+).\]

\begin{lemma} \label{L:radial}
 $M_k^1(\b R^N)$ is a commutative semigroup with respect to $*_k$ with 
neutral element $\delta_0$. 
The mapping $\,p: (M_k^1(\b R^N), *_k)\,\to\, (M^1(\b R_+), \circ_\lambda), \,\,\mu\mapsto 
        p(\mu)$ establishes
 an isometric isomorphism of semigroups as well as a homeomorphism with
respect to the weak topologies on both spaces.
\end{lemma}

\begin{proof} 
It is obvious that $p$ is an isometric bijection from $M_k^1(\b R^N)$ onto
$M^1(\b R_+)$. Moreover,  by 
L\'evy's continuity theorem for the Dunkl transform and the 
Hankel transform respectively,
we deduce from \eqref{(5.110)} that $p$ is a homeomorphism w.r.t. the
weak topologies on the stated spaces. Commutativity of $*_k$ 
and the statement on the neutral element are clear. It remains to prove that
$M_k^1(\b R^N)$ is closed with respect to $*_k$, and that $p$ is multiplicative. For this, 
let
 $\mu, \nu\in M_k^1(\b R^N).$ Then by \eqref{(5.101)} 
and \eqref{(5.110)},
\begin{align}
(\mu*_k\nu)^{\wedge k}(\xi) \,=\,& \widehat \mu^{\,k}(\xi)\widehat \nu^{\,k}(\xi) \,=\,
  \mathcal H^\lambda(p(\mu))(|\xi|)\mathcal H^\lambda(p(\nu))(|\xi|) \notag\\
  =\,&\mathcal H^\lambda(p(\mu)\circ_\lambda p(\nu))(|\xi|)\notag
\end{align}
Let $\tau\in M_k^1(\b R^N)$ with $\,p(\tau) = p(\mu)\circ_\lambda p(\nu).$ 
Then the above calculation shows that $\, (\mu*_k \,\nu)^{\wedge k}\,=\, \widehat \tau^{\, k}.$
The injectivity of the Dunkl transform  implies  that
 $\,\mu*_k\nu = \tau \in M_k^1(\b R^N),$ as well as
\[ p(\mu*_k\nu) = p(\tau) = p(\mu)\circ_\lambda p(\nu).\]
This finishes the proof.
\end{proof}

\noindent
Recall that for two Markov kernels $P,  Q$ on $\b R^N$, the composition $P\circ Q$ is defined by
\[ P\circ Q(x,A) = \int_{\b R^N} Q(z,A) P(x,dz).\]
If $P$ and $Q$ are $k$-invariant, then  $P\circ Q$ is again $k$-invariant   with 
\[ \bigl((P\circ Q)(x,\,.\,)\bigr)^{\wedge k} (\xi) =\,P(0,\cdot\,)^{\wedge k}(\xi)\cdot 
Q(0,\cdot\,)^{\wedge k}(\xi)\cdot E_k(-ix,\xi) \quad (x,\xi\in \b R^N),\]
c.f. \cite{RV}. A family $(P_t)_{t\geq 0}$ of $k$-invariant Markov kernels on $\b R^N$ is called a 
semigroup of $k$-invariant Markov kernels, if the following statements hold:
\begin{enumerate}
\item[\rm{(i)}] The kernels $(P_t)_{t\geq 0}$ form a semigroup, i.e. $\,P_s\circ P_t = P_{s+t}$ for $s,t\geq 0.$ 
\item[\rm{(ii)}] The mapping $\b R_+ \to M^1(\b R^N), \,t\mapsto P_t(0,.\,)$ is weakly continuous.
\end{enumerate} 

The above results enable us to construct   
semigroups of $k$-invariant Markov kernels from
 convolution semigroups  with respect to the hypergroup structure $\circ_\lambda$ on 
$\b R_+$. 

\begin{definition}
A family $(\sigma_t)_{t\geq 0} \subset M^1(\b R_+)$ of 
probability measures on $\b R_+$ is called a 
convolution semigroup on $\b R_+$ with respect to $\circ_\lambda$, if 
$\sigma_s\circ_\lambda\sigma_t = \sigma_{s+t}$ for all $s,t\geq 0$ with 
$\sigma_0 = \delta_0$, and if the mapping $\,\b R_+ \to M^1(\b R_+), \,
t\mapsto \sigma_t$ is weakly continuous.
\end{definition}

\begin{theorem}\label{T:semigroup}
 Suppose that $(\sigma_t)_{t\geq 0}$ is a convolution semigroup on 
$\b R_+$ w.r.t. $\circ_\lambda$, and  define $(\mu_t)_{t\geq 0}\subset M^1_k(\b R^N)$  
by $\,\mu_t:= p^{-1}(\sigma_t)\in M_k^1(\b R^N)$. 
Then $\, \mu_s *_k\mu_t = \mu_{s+t}$ for all $s,t\geq 0$, and 
$\, P_t(x,A):= \delta_x\circ_k\mu_t \,$
defines a semigroup of $k$-invariant Markov kernels on $\b R^N$. 
\end{theorem}

\begin{proof}
Let $s,t\geq 0.$ Then by Lemma \ref{L:radial},
\[p(\mu_s*_k\mu_t)\,=\,\nu_s\circ_\lambda\nu_t  \,=\, \nu_{s+t} \,=\, p(\mu_{s+t}).\]
As $p$ is injective on $M_k^1(\b R^N),$ this proves that
$\, \mu_s *_k\mu_t = \mu_{s+t}$. Further, the weak continuity of
the mapping $t\mapsto \mu_t$ on $\b R_+$ is clear from the continuity of $p^{-1}$. 
According to Theorem \ref{T:Markov2},
each $P_t$ is a $k$-invariant
Markov kernel with $P_t(0,\,.\,) = \mu_t$,  and 
\begin{align} (P_s\circ P_t)(x,.\,)^{\wedge k}& (\xi)\,=\, E_k(-ix,\xi) P_s(0,.\,)^{\wedge k}(\xi)\,
P_t(0,.\,)^{\wedge k}(\xi)\notag\\
=\,& E_k(-ix,\xi) \widehat \mu_s^{\,k}(\xi)\widehat \mu_t^{\,k}(\xi)\,=\,
 E_k(-ix,\xi)  \widehat \mu_{s+t}^{\,k}(\xi)\,=\, P_{s+t}(x,.\,)^{\wedge k}(\xi).\notag
\end{align}
This implies the semigroup property of $(P_t)_{t\geq 0}$.
\end{proof}

\begin{examples}\label{E:ex}
All the examples for semigroups of $k$-invariant 
Markov kernels discussed in \cite{RV} 
fit into the concept of Theorem \ref{T:semigroup}: 

\noindent
(1) ($k$-Gaussian semigroups.) The $k$-Gaussian semigroup 
$(P_t^\Gamma)_{t\geq 0}$ is defined by
\[ P_t^\Gamma(x,A) = \, \int_A \Gamma_k(t,x,y)w_k(y)dy\,,\quad(t>0)\]
with the Dunkl-type heat kernel $\Gamma_k$ (c.f. \eqref{(3.102)}).
It is $k$-invariant and defines a  Feller process on $\b R^N$, 
the Dunkl-type Brownian motion. 
According to \eqref{(3.101)} and Lemma \ref{L:Translation}, we may write
\[ \Gamma_k(t,x,y) = F_k(t,-x*_ky) \]
with
\[F_k(t,x):= \Gamma_k(t,x,0)\,=\,   
           \frac{1}{(2t)^{\lambda +1}c_k}e^{-|x|^2/4t}\,.\]
For $t >0,$ put
\[ d\mu_t(x):= F_k(t,x) w_k(x)dx \in M_k^1(\b R^N).\]
Then $\sigma_t:= p(\mu_t)\in M_\gamma^1(\b R_+)$ is given by the Rayleigh 
distribution
\[ d\sigma_t(r) = \frac{d_k}{(2t)^{\lambda +1}c_k}r^{2\lambda +1}
 e^{-r^2/4t}, \, \, t>0.\]
It is well-known (see e.g. 7.3.18 of \cite{BH}) that $(\nu_t)_{t\geq 0}$ (with $\nu_0:=0$) is a convolution semigroup on $\b R_+$ w.r.t. 
 $\circ_\lambda$. It defines
a Bessel process of index $\lambda$. In order to see that 
$(P_t^\Gamma)_{t\geq 0}$ is associated with $(\sigma_t)_{t\geq 0}$ as in Theorem \ref{T:semigroup}, it remains to verify that
\[ \delta_x *_k\mu_t \,=\, P_t^\Gamma(x,\,.\,).\]
But if $f\in \mathscr S(\b R^N),$ then by use of Lemma \ref{L:Translation}
\[ \delta_x *_k\mu_t(f) = \, \int_{\b R^N} F_k(t,y) f(x*_ky) w_k(y)dy\,=\,
   \int_{\b R^N} \Gamma_k(t,x,y) f(y)w_k(y)dy.\]
This yields the claimed identity. 

\smallskip
\noindent
(2) (Subordination) 
If $(P_t)_{t\geq 0}$ is a semigroup of $k$-invariant Markov kernels and $(\rho_t)_{t\geq 0}\subset M^1(\b R)$ is a 
convolution semigroup of
probability measures on the group $(\b R, +)$ (in the sense of \cite{BF}) which
is supported by $\b R_+$, then  a
 subordinated semigroup $(\widetilde P_t)_{t\geq 0}$ 
of $k$-invariant Markov kernels is defined by
\[\widetilde P_t(x,A) = \, \int_0^\infty P_s(x,A) d\rho_t(s).\]
Suppose now that  $(P_t)_{t\geq 0}$ is 
associated with a $\circ_\lambda$-
convolution semigroup $(\sigma_t)_{t\geq 0}$ on $\b R_+$
according to Theorem \ref{T:semigroup}. 
Then it is immediate that $(\widetilde P_t)_{t\geq 0}$ is associated
in the same way with the  $\circ_\lambda$-
convolution semigroup $(\widetilde \sigma_t)_{t\geq 0}$ defined by
\[ \widetilde \sigma_t := \int_0^\infty \sigma_s\, d\rho_t(s).\]
This example in particular includes the $k$-Cauchy kernels in
\cite{RV}.

\end{examples}

\section{Appendix: Generalized Bessel functions as spherical functions}\label{euclidean}

For crystallographic reflection groups  
and certain discrete sets of half-integer multiplicity functions, 
generalized Bessel functions  
have an interpretation as the spherical functions of a 
Euclidean-type symmetric space.  This in particular implies a positive product formula
as stated in \eqref{(1.1f)}.

To become more precise, let us first recall the context from \cite{O} and 
\cite{dJ2}:
Suppose that 
$\mathcal G$ is a connected, non-compact semisimple Lie group 
with finite center. Choose a  maximal compact subgroup $\mathcal K$ of 
$\mathcal G$, and let  $\mathcal K\ltimes \frak p$ be the  Cartan motion group
 associated with the Cartan decomposition $\frak g = \frak k + \frak p$ of
the Lie algebra $\frak g$; here 
$\mathcal K$ acts on $\frak p$ via the adjoint representation. 
Choose a maximal abelian subspace $\frak a$ of $\frak p$, and denote by  
by $\Sigma\subset \frak a^*$ the roots of $\frak g$ 
with respect to $\frak a$. 
We consider $\Sigma$ as a subset of $\frak a$, identifying 
$\frak a$ with its dual 
via the Killing form $B$, and denote by $G$ the Weyl group of $\Sigma$, 
acting on the Euclidean space $(\frak a, B) \cong (\b R^N, \langle\,.\,,\,.\,\rangle)$.
Let us consider the spherical functions of the Gelfand pair 
 $(\mathcal K \ltimes \frak p, \mathcal K)$  as $\mathcal K$-invariant 
functions on $(\mathcal K\ltimes \frak p)/\mathcal K \,\cong \, \frak p$.  
According to Chap. IV of  \cite{Hel}, they are of the form $\psi_\lambda\,,\,\lambda\in 
\frak a_{\b C}$, where 
$\psi_\lambda\in C^\infty(\frak p)$
 is  characterized as the unique $\mathcal K$-invariant
 solution of the 
joint eigenvalue
problem
\[ \partial(p)\psi\,=\, p(\lambda) \psi\,\, \text{ for all }p\in I(\frak p),\quad \psi(0) =1;\]
here  $I(\frak p)$ is the 
space of $\mathcal K$-invariants in the symmetric algebra over $\frak p$ and $\partial(p)$ 
is the constant coefficient differential operator associated with $p$. Further, $\psi_\lambda = \psi_\mu$ 
iff $\mu = g\lambda$ for some $g\in G$. 
By taking radial parts of constant coefficient differential operators 
on $\frak p$, it can now be shown that the  restriction
$\psi_\lambda\vert_\frak a$ (which is $G$-invariant and determines $\psi_\lambda$ uniquely) 
coincides with
the a generalized Bessel function $J_k(\,.\,,\lambda)$ 
associated with $G$;   the multiplicity function $k$ is given by 
$ k(\alpha) = \frac{1}{4} \sum_{\beta\in \b R\alpha\cap\Sigma} m_\beta$, 
with  $m_\alpha$ the multiplicity of $\alpha \in \Sigma$.  A detailed proof of this 
can be found in 
\cite{dJ2}.

The convolution of 
$\mathcal K$-biinvariant functions on the Cartan
motion group can be interpreted as the convolution
of a commutative hypergroup structure on the double coset space  
$(\mathcal K\ltimes \frak p)/\!/\mathcal K\,\cong \mathcal K\setminus\! \frak p$ which is defined
 by 
\[ \delta_{\mathcal K\cdot x} * \delta_{\mathcal K\cdot y}\,=\,
\int_{\mathcal K} \delta_{\mathcal K\cdot(x+k\cdot y)}\,dk, \quad x,y\in \frak p.\]
(See \cite{J}, Chap.~8  for information about double coset and orbit hypergroups.)
If  $\frak a_+$ is an arbitrary fixed Weyl chamber of $\frak a$, then 
each orbit $\mathcal K\cdot x$ contains a unique $x_+\in \overline{\frak a_+}$, and 
 the mapping
$\, \mathcal K\setminus\! \frak p   \mapsto \overline{\frak a_+}, \,\, \mathcal K\cdot x \to x_+\,$
is a homeomorphism; see \cite{Hel}, Prop. I.5.18. Thus, we obtain an induced hypergroup structure on the closed chamber $\overline{\frak a_+}\,$ which 
is determined by the convolution \begin{equation}\label{(7.1)}
 \delta_x * \delta_y \,=\, \int_{\mathcal K} \delta_{(x+k\cdot y)_+}dk\,, \,\,\, x,y\in 
           \overline{\frak a_+}.
\end{equation}
In particular, the $\delta_x*\delta_y$ are compactly supported probability measures on 
$\overline{\frak a_+}$.
The identity element of this hypergroup is $0$, and the involution is given by $(x_+)^- = ((-x)_+)^{-}.$
From \eqref{(7.1)} together with the characterization of the spherical functions $\psi_\lambda$
by means of their product formula, we see that the set
\[ \{\psi\in C(\overline{\frak a_+}\,): \delta_x*\delta_y(\psi) = \psi(x)\psi(y)\}\]
 of continuous, multiplicative 
functions on the hypergroup $\overline{\frak a_+}$ consists exactly of the 
functions $\psi_\lambda\,, \,\lambda\in \frak a_{\b C}$. With the identification 
 $(\frak a, B) \cong (\b R^N, \langle\,.\,,\,.\,\rangle), \,
\frak a_+ \cong C$ and $J_k$ as defined above, we 
in particular  obtain
\[ J_k(x,\lambda) J_k(y,\lambda) = \int_{\overline{C}} J_k(\xi,\lambda) 
 d(\delta_x*\delta_y)(\xi)\quad\, \text{for all }\,x,y\in \overline C, \, \lambda\in \b C^N.\]
Suppose now that $\mathcal G$ is complex. Then $k(\alpha) =1$ for all $\alpha\in \b R$, and 
the associated Bessel function $J_{\bf 1}(\,.\,,\lambda)$
is given by 
\[ J_{\bf 1}(x,\lambda)\,=\, c\sum_{g\in G} \frac{\text{det}\,(g)}{\pi(x)\pi(\lambda)}\,
 e^{\langle x,g\lambda\rangle}\]
with some constant $c\in \b C$; here $\pi$ denotes
 the fundamental alternating polynomial
\[ \pi(x)= \prod_{\alpha\in R_+} \langle\alpha,x\rangle\,,\]
c.f. \cite{D3}, Prop. 1.4. (This result  remains obviously 
true if not a single $G$-orbit in $R$, but any $G$-invariant subset, in particular $R$ itself, is considered).  
The Laplace representation \eqref{(2.77)} for $E_k$ implies
that for all $x,\lambda\in \b R^N$, 
\[ J_{\bf 1}(x,\lambda)\,=\, \int_{\b R^N} e^{\langle\xi,\lambda\rangle} 
d\widetilde \mu_x(\xi)\]
with the $G$-invariant probability measure
\[ \widetilde\mu_x:= \frac{1}{|G|}\sum_{g\in G} \mu_{gx}^{\bf 1}.\]
For $y\in \b R^N$, denote further 
by $\nu_{x,y}\in M^1(\b R^N)$ the image measure of $\widetilde\mu_x$
under the translation $\xi\mapsto\xi+y$.
 Then 
\begin{align}
 J_{\bf 1}(x,\lambda) J_{\bf 1}&(y,\lambda)\,=\,c\int_{\b R^N} \sum_{g\in G} 
\frac{\text{det}\,(g)}{\pi(y)\pi(\lambda)}
 e^{\langle y+\xi, g\lambda\rangle} d\widetilde\mu_x(\xi)\notag\\ 
=\,& c\int_{\b R^N} \sum_{g\in G} \frac{\pi(\xi)}{\pi(y)}\frac{\text{det}\,(g)}{\pi(\xi)\pi(\lambda)} e^{\langle\xi, g\lambda\rangle}d\nu_{x,y}(\xi)\,=\,
\int_{\b R^N}J_{\bf 1}(\xi,\lambda)\frac{\pi(\xi)}{\pi(y)} d\nu_{x,y}(\xi).\notag
\end{align}
This shows that in the present case, the hypergroup convolution on $\overline C$ is given by
$\, \delta_x*\delta_y = \pi(y)^{-1}\rho_{x,y},$ where $\rho_{x,y}$ is the image measure
of $\pi\nu_{x,y}$ under the canonical projection $\b R^N\to G\!\setminus\!\b R^N \cong \overline C$. 
Up to a common multiplictive factor, the spherical functions $\psi_\lambda$ of 
 $(\mathcal K \ltimes \frak p, \mathcal K) $  
can be identified with the spherical functions $\Phi_\lambda$ of the symmetric space
$\mathcal G/\mathcal K$, c.f. \cite{Hel}, Prop. IV.4.10. So the $\Phi_\lambda$ 
 satisfy a slightly modified
positive product formula, which  has recently been investigated  in \cite{GS}.
The representing measures are absolutely continuous with respect to Lebesgue measure in the generic case, and  the authors obtain detailed information on their supports.

\bibliographystyle{amsplain}

\end{document}